\documentclass[]{article}

\usepackage{graphicx}
\usepackage{amsbsy,bm,relsize}
\usepackage[affil-it]{authblk}
\usepackage{url}

\begin{document}
\title{Unavoidable golden ratio}

\author{Dorota Jacak}

\affil{{\small Institute of Mathematics and Computer Science, Wroclaw University of Technology, Wyb. Wyspia{\'n}skiego 27, 50-370 Wroc{\l}aw, Poland\\dorota.jacak@pwr.wroc.pl}}

\maketitle

\begin{abstract}
It is demonstrated that iterative repeating of some simple geometric construction leads unavoidably in the limit to the golden ratio. The procedure appears to be quickly convergent regardless of a ratio of initial elements sizes. This could explain a widespread occurrence of the golden ratio in various constructions, including puzzling proportions encountered in structures of ancient cultures (e.g., in the Great Pyramid of Giza), as well as in the natural growth and self-organization processes in organic and inorganic matter.
\end{abstract}

\paragraph{Introduction}
It is well known, that the golden ratio frequently occurs in various constructions, including puzzling proportions encountered in structures of ancient cultures, as well as in the natural growth processes and self-organization in organic and inorganic matter \cite{www, quant}. The cross-section of the Great Pyramid of Giza (pyramid of Cheops) suggests that this structure is based on the so-called Kepler triangle, i.e., the right triangle with the sides, $a$, $a\sqrt{\phi}$, $a\phi$, $a>0$, where $\phi=1.618033989\ldots$ is the golden ratio, and this proportion is repeated in the pyramid with the surprisingly high precision. Although the occurrence of the golden ratio, or rather an almost-ideal golden ratio is widespread in the nature, the reason of this fact is still unclear. In the present letter we demonstrate that iterative repeating of some special construction steps always leads in the limit to the golden ratio, regardless of a size ratio of elements assumed as the initial conditions. Proving of a simple theorem, on which few simple constructions can be based, gives a possible explanation of so vast manifestation of the golden ratio. The remarkable feature of suggested constructions is their very quick convergence, i.e., only after a very few steps they arrive at almost precise approximation of the golden ratio regardless to initial geometric proportions, even far from the golden ratio. As the constructions require only few repeating tasks of simple measurements, they are feasible with usage of even primitive tools (like rope, etc.) and one could guess that those methods might have been practiced by the ancient Egyptians. The simplicity and quick convergence of the procedure suggest also that self-organised processes and natural growth would employ this scheme, which results in the abundance of the golden ratio patterns in corresponding step-by-step built structures of many natural systems.

\paragraph{Central theorem}

For the Fibonacci sequence, $\left(\phi(k)\right)=\left\{1,1,2,3,5,8,\ldots \right\}$, $\phi(k) + \phi(k+1) = \phi(k+2)$, one can write \cite{vajda}:
\begin{equation}
\label{eq1}
\phi = \lim_{k\rightarrow \infty}\frac{\phi(k+1)}{\phi(k)},
\end{equation}
where $\phi=1.618033989\ldots$ is the golden ratio. This series is quickly convergent -- just for $k \sim 5$ it gives almost the limiting value $\phi$ (e.g., $\frac{\phi(6)}{\phi(5)}=1.6$).\\

{\it Theorem:} Let $(p_{k})$ be the series  of numbers defined by the following recurrent formula:
\begin{equation}
\label{eq2}
p_0=a, p_1=\lambda a, a,\lambda>0, \displaystyle\mathop{\mathlarger{\mathlarger{\mathlarger{\forall}}}}_{k=2,3,4,\ldots}p_k=\sqrt{p_{k-1}^2+p_{k-2}^2},
\end{equation}
then
\begin{equation}
\lim_{k\rightarrow \infty}\frac{p_{k+1}}{p_{k}}=\sqrt{\phi}.
\end{equation}\\

{\it Proof:} By the induction (assuming $\phi(0) = 0$), one can show that 
$\displaystyle\mathop{\mathlarger{\mathlarger{\mathlarger{\forall}}}}_{k=1,2,3,\ldots}p_k = a\sqrt{\phi(k-1)+\lambda^2\phi(k)}$, which ends the proof by virtue of eq. (\ref{eq1}).

\paragraph{Constructions}

\begin{figure}[!h]
\centering
\scalebox{1.0}{\includegraphics{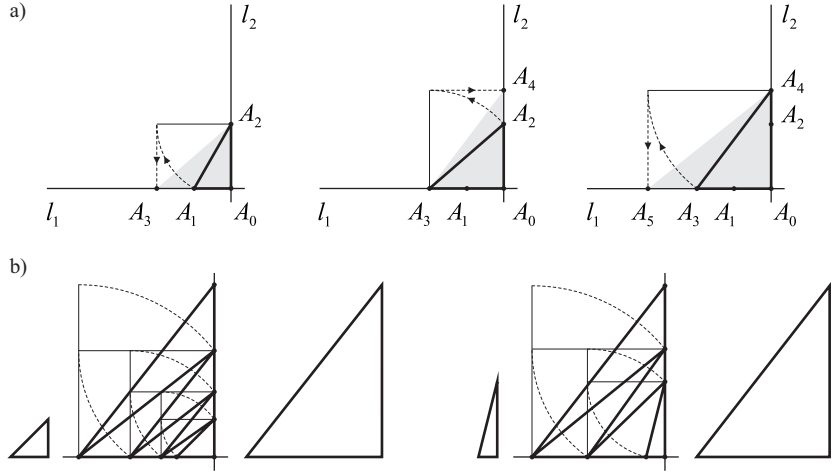}}
\caption{\label{fig1} {\footnotesize a) Three initial steps of the construction of an almost-ideal Kepler triangle (in the odd steps the length of the base is increased to the length of the hypotenuse of the triangle from the former step and the height remains unchanged; in the even steps the length of the height is increased to the length of the hypotenuse of the triangle from the former step and the base remains unchanged) -- the grey area is the newly constructed triangle in each step, the bold line draws the triangle from a former step (in the first step this is the initial triangle);  the hypotenuse at each step marks point on one of perpendicular lines creating the series $A_1A_2 = A_0A_3, A_2A_3 = A_0A_4, A_3A_4 = A_0A_5$; b) the construction leads to the Kepler triangle regardless of the form of initial triangles}}
\end{figure}

Using the above theorem one can easily construct an almost-ideal Kepler triangle. The procedure resolves to creation of the series of rectangular triangles, one in each consecutive step, in such a way that: in the odd steps the length of the base is increased to the length of the hypotenuse of the triangle from the former step and the height remains unchanged; in the even steps the length of the height is increased to the length of the hypotenuse of the triangle from the former step and the base remains unchanged (cf. Figure \ref{fig1}). The more formal construction is as follows. Let two perpendicular lines, $l_0$ and $l_1$, intersect in the point $A_0$. Let us  mark two of arbitrary length  segments along these lines: $\overline{A_0A_1},A_0A_1=a, A_1\in l_1, \overline{A_0A_2},A_0A_2=\lambda a, A_2\in l_2, a, \lambda > 0$. Next, let us mark points $A_{2k+1}\in l_1, A_{2k}\in l_2, k\geq 1$, in such way that $A_0A_{2k+1}=A_{2k-1}A_{2k}, A_0A_{2k}=A_{2k-2}A_{2k-1}$, as it is depicted in Figure \ref{fig1}.

The sequence $\left(A_0A_k\right)$ satisfies the condition (\ref{eq2}) and thus, $\lim_{k \rightarrow \infty}\frac{A_0A_{k+1}}{A_0A_k} = \sqrt{\phi}$. Each of the triangles $A_{2k-1}A_0A_{2k}$ and $A_{2k}A_0A_{2k+1}$ is an almost-ideal Kepler triangle already just for $k \geq 4$.

If one assumes $A_0A_1 = 3$ and $A_0A_2 = 4$ (in a correspondence to a common measurement tool in ancient Egypt in the form of a rope with knots apart by 3, 4 and 5 length units, allowing to create a 3-4-5 right triangle  known as the Egyptian triangle) one could obtain an almost-ideal Kepler triangle in a few consecutive  steps, because the ratio of just 8-th to 7-th hypotenuse (obtained in the above described construction) gives a value $\approx 1.272$ ($\sqrt{\phi} = 1.272019649\ldots$).

\begin{figure}[!h]
\centering
\scalebox{1.0}{\includegraphics{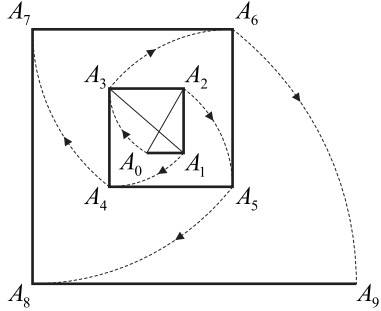}}
\caption{\label{fig2} {\footnotesize The construction of an almost-ideal golden polygonal chain}}
\end{figure}

In the similar way as an almost-ideal Kepler triangle is constructed, one can construct a polygonal chain (cf. Figure \ref{fig2}). It is constructed by consecutive connecting of points $A_0,A_1,A_2,\ldots$ on the plain in such manner that $\overline{A_0A_1}\perp\overline{A_1A_2}$, $A_0A_1=a$, $A_1A_2=a\lambda$, $a,\lambda>0$, $\overline{A_{k-1}A_k}\perp\overline{A_kA_{k+1}}$ and $A_kA_{k+1}=A_{k-2}A_k=\sqrt{\left(A_{k-2}A_{k-1}\right)^2+\left(A_{k-1}A_k\right)^2}$, 
$k=2,3,4,...$. Vertices of polygonal chain are arranged on the plain in counterclockwise ordering. The sequence $\left(A_kA_{k+1}\right)$ satisfies the condition (\ref{eq2}), thus $\lim_{k \rightarrow \infty}\frac{A_kA_{k+1}}{A_{k-1}A_k} = \sqrt{\phi}$ and $\lim_{k \rightarrow \infty}\frac{A_kA_{k+1}}{A_{k-2}A_{k-1}} = \phi$, so the triangle $A_{k-1}A_kA_{k+1}$ is an almost-ideal Kepler triangle already just for $k\geq 8$.

If $\lambda = \sqrt{\phi}$, then such a polygonal chain can be called the golden polygonal chain. For the golden polygonal chain any triangles $A_{k-1}A_kA_{k+1}$ and $A_{k-1}SA_k, k=1,2,3,\ldots$, are the Kepler triangles, where $S = l_1\cap l_2$, line $l_1:A_0,A_2\in l_1$, line $l_2:A_1,A_3\in l_2$. It is easy to notice that $l_1$ and $l_2$ lines are perpendicular, $A_{2k}\in l_1$, $A_{2k+1}\in l_2$ and $S$ divides the segments $\overline{A_{2k}A_{2k+2}}$ and $\overline{A_{2k+1}A_{2k+3}}$, $k=0,1,2,...$, according to the golden ratio.

\begin{figure}[!h]
\centering
\scalebox{1.0}{\includegraphics{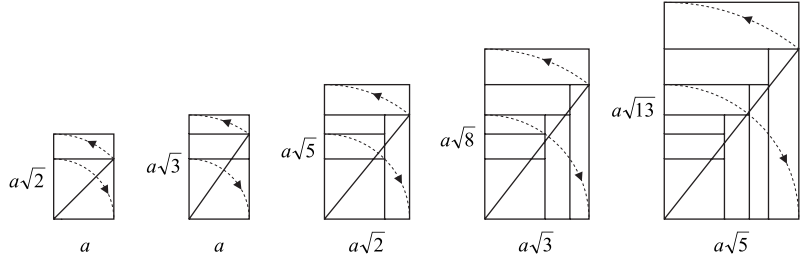}}
\caption{\label{fig3} {\footnotesize The construction of an almost-ideal golden rectangle. Could this have been used in elements of the Alhambra structures?}}
\end{figure}

\begin{figure}[!h]
\centering
\scalebox{1.0}{\includegraphics{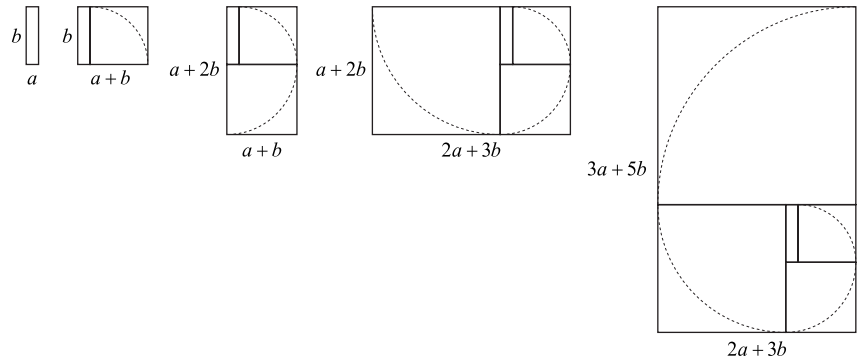}}
\caption{\label{fig4} {\footnotesize The almost-ideal golden rectangle construction steps -- clockwise consecutive appending of a square (of side equal to a longer rectangle side) to a longer side of a rectangle. After repeating the construction steps one will obtain in the limit an ideal golden rectangle and simultaneously the dashed curve will form the ideal golden spiral (as for a finite number of steps it will be an almost-ideal golden spiral).}}
\end{figure}

Besides the construction of an almost-ideal Kepler triangle one can also propose similarly easy constructions of an almost-ideal golden rectangle (cf. Figures \ref{fig3} and \ref{fig4}). In the construction illustrated in Figure \ref{fig3} after $k$ step one obtains a rectangle with sides $a\sqrt{\phi(k)}$ and $a\sqrt{\phi(k+2)}$, which is an almost-ideal golden rectangle, ideal in the limit, because, $\lim_{k \rightarrow \infty}\frac{\sqrt{\phi(k+2)}}{\sqrt{\phi(k)}} = \phi $. Again the series is quickly convergent. 

In the $k$-th step of construction illustrated in Figure \ref{fig4} (if $b \geq a$) one obtains the rectangle with the sides $a\phi(k-1)+b\phi(k)$, $a\phi(k)+b\phi(k+1)$ and the following equalities are satisfied, $\lim_{k \rightarrow \infty}\frac{a\phi(k)+b\phi(k+1)}{a\phi(k-1)+b\phi(k)} = \phi$, thus that rectangle is an almost-ideal golden rectangle.

The remarkable property is that all of the above presented constructions are effective within the first few steps regardless of a size and mutual ratio of elements assumed as  initial ones. This might be responsible for widespread manifestations of the golden ratio in the natural self-organised processes, especially in the growth process schematically presented in Figures \ref{fig2} and \ref{fig4}.

\bibliography{dj}
\bibliographystyle{unsrt}

\end{document}